\documentclass[reqno,12pt]{amsart}

\numberwithin{equation}{section}

\usepackage{verbatim} 
\usepackage{xspace} 
\usepackage{amssymb}
\usepackage{amsmath}
\usepackage{amsthm} %
\usepackage[latin1]{inputenc}
\usepackage{graphicx}
\usepackage{amscd,amssymb,euscript,graphics}
\makeindex

\newcommand{\C}{\ensuremath{\mathbb{C}}}
\newcommand{\F}{\ensuremath{\mathbb{F}}}
\newcommand{\N}{\ensuremath{\mathbb{N}}}

\newcommand{\norm}[1]{\ensuremath{\|#1\|}}

\newcommand{\bc}{\begin{center}}
\newcommand{\ec}{\end{center}}
\newcommand{\be}{\begin{enumerate}}
\newcommand{\ee}{\end{enumerate}}
\newcommand{\bi}{\begin{itemize}}
\newcommand{\ei}{\end{itemize}}
\newcommand{\bd}{\begin{description}}
\newcommand{\ed}{\end{description}}
\newcommand{\beq}{\begin{equation}}
\newcommand{\eeq}{\end{equation}}
\newcommand{\beqa}{\begin{eqnarray}}
\newcommand{\eeqa}{\end{eqnarray}}
\newcommand{\bfr}{\begin{flushright}}
\newcommand{\efr}{\end{flushright}}
\newcommand{\bfl}{\begin{flushleft}}
\newcommand{\efl}{\end{flushleft}}

\newcommand\Act{{\widetilde{\Ac}}}
\newcommand\Bt{{\widetilde{B}}}
\newcommand\Et{{\widetilde E}}
\newcommand\ImPt{{\mathrm{Im}\;}}
\newcommand\lambdat{{\tilde\lambda}}
\newcommand\RePt{{\mathrm{Re}\;}}
\newcommand\vol{{\operatorname{vol}}}
\newcommand\Yt{{\widetilde{Y}}}
\newcommand\Zt{{\widetilde{Z}}}

\newtheorem{teo}{Theorem}[section]

\newtheorem{lem}[teo]{Lemma}

\theoremstyle{remark}

\theoremstyle{definition}

\newtheorem{df}[teo]{Definition}
\newtheorem{obs}[teo]{Remark}

\newcommand\Ac{{\mathcal{A}}}

\newcommand{\bb}[1]{{\mathbb{#1}}}
\newcommand\Bc{{\mathcal{B}}}

\newcommand{\cl}[1]{{\mathcal{#1}}}
\newcommand\Cpx{\bb{C}}

\newcommand\Dc{{\mathcal{D}}}

\newcommand\DT{\operatorname{DT}}
\newcommand\eps{\epsilon}

\newcommand\Gc{{\mathcal{G}}}
\newcommand\Hc{{\mathcal{H}}}

\newcommand\Mcal{\cl M}

\newcommand\Nc{\cl N}

\newcommand\tr{{\operatorname{tr}}}

\newcommand\Wc{{\mathcal{W}}}

\newcount\theTime
\newcount\theHour
\newcount\theMinute
\newcount\theMinuteTens
\newcount\theScratch
\theTime=\number\time
\theHour=\theTime
\divide\theHour by 60
\theScratch=\theHour
\multiply\theScratch by 60
\theMinute=\theTime
\advance\theMinute by -\theScratch
\theMinuteTens=\theMinute
\divide\theMinuteTens by 10
\theScratch=\theMinuteTens
\multiply\theScratch by 10
\advance\theMinute by -\theScratch

\def\today{{\number\day\space
 \ifcase\month\or
  January\or February\or March\or April\or May\or June\or
  July\or August\or September\or October\or November\or December\fi
 \space\number\year}}

\begin{document}

\title[Free entropy dimension]
{The free entropy dimension of some $L^{\infty}[0,1]$-circular operators}

\author[Ken Dykema and Gabriel Tucci]
{Kenneth J.\ Dykema*, Gabriel H. Tucci}

\address{Department of Mathematics, Texas A\&M University,
College Station, TX 77843-3368, USA}
\email{kdykema@math.tamu.edu}
\email{gtucci@math.tamu.edu}

\thanks{*Research supported in part by NSF grant DMS--0300336.}


\begin{abstract}
We find the microstates free entropy dimension of a large class of
$L^\infty[0,1]$--circular operators, in the presence of a generator of the diagonal subalgebra.
\end{abstract}

\maketitle

\section{Introduction}

Let $\Mcal$ be a von Neumann algebra with a specified normal faithful tracial state $\tau$.
The free entropy dimension
\begin{equation}\label{eq:delta0}
\delta_0(X_1,\ldots,X_n)
\end{equation}
for $X_1,\ldots,X_n\in\Mcal$,
was introduced by Voiculescu~\cite{V94}, \cite{V}, see also~\cite{V02}.
This quantity is sometimes called
the microstates free entropy dimension to distinquish it from
another version introduced by Voiculescu and because its definition utilizes
matricial microstates for the operators $X_1,\ldots,X_n$.
It is an open problem whether the quantity~\eqref{eq:delta0} is an invariant
of the von Neumann algebra generated by $X_1,\ldots,X_n$,
and it is of interest to find the free entropy dimension of various operators.
See, for example~\cite{V}, \cite{V99} \cite{GS02}, \cite{J03}, \cite{DT},
\cite{J-hyp}, \cite{J-1bd}, \cite{J-subf}, \cite{JS} for some such results.

\begin{figure}[Ht]
\begin{center}
\includegraphics[width=4cm]{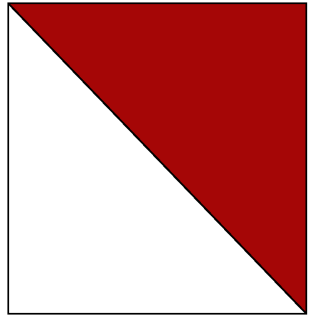}
\caption{The upper trangle, representing the quasinilpotent DT-operator $T$}\label{figure1}
\end{center}
\end{figure}

In~\cite{DT}, Dykema, Jung and Shlyakhtenko computed $\delta_0(T)=2$
for the quasinilpotent DT-operator $T$.
This operator was introduced by Dykema and Haagerup in~\cite{DT0}.
It can
be realized as a limit in $*-$moments of strictly upper-triangular
random matrices with i.i.d.\ complex Gaussian entries above the
diagonal.
Alternatively, as was seen in \cite{DT0}, $T$ can be
obtained in the free group factor $L(\F_2)$
from a semicircular element $X$ and a free copy of
$L^{\infty}([0,1])$ by using projections from the latter to cut
out the upper triangular part of $X$.
(Note that $X$ may be
replaced by a circular element $Z$ for this procedure.)
Then we
can visualize $T$ as in Figure~\ref{figure1},
where the shaded
region has weight 1, the unshaded region has weight 0, and these
weights are used to multiply entries of a Gaussian random matrix.
It was proved in
\cite{DT1} that the von Neumann algebra generated by $T$ contains
all of $L^\infty([0,1])$, and is, thus,
the free group factor $L(\F_2)$.

In this paper we consider more general operators than $T$, defined also
as limits of random matrices or,
equivalently, in the approach was taken in \cite{B-circ}, by
cutting a circular operator $Z$ using projections in a $*$--free copy
of $L^{\infty}([0,1])$.
The class of operators considered there
consisted of those $L^{\infty}([0,1])$-circular operators described as
follows.
Let $\eta$ be an absolutely continuous measure with respect
to Lebesgue measure on $[0,1]^2$ with Radon--Nikodym derivative
$H\in L^{1}([0,1]^2)$ and assume the push--forward measures $\pi_{i*}\eta$ under
the coordinate projections $\pi_{1},\pi_{2}:[0,1]^2\to [0,1]$ are
absolutely continuous with respect to Lebesgue measure and have
essentially bounded Radon--Nikodym derivatives.
For each such measure
$\eta$ with the associated function $H\in L^{1}([0,1]^2)$
we have the operator
$Z_{H}$ described in~\cite{B-circ};
(however, this operator was denoted $z_{\eta}$ in~\cite{B-circ}).
When $\eta$ is Lebesgue measure on $[0,1]^2$, then $H=1$ and
$Z_{H}$ is the usual circular operator. When $\eta$ is the
restriction of Lebesgue measure to the upper triangle pictured in
Figure~\ref{figure1}, then $H$ is the characteristic function of this triangle
and $Z_{H}$ is the
quasinilpotent DT-operator $T$.

\begin{figure}[Ht]
\begin{center}
\includegraphics[width=4cm]{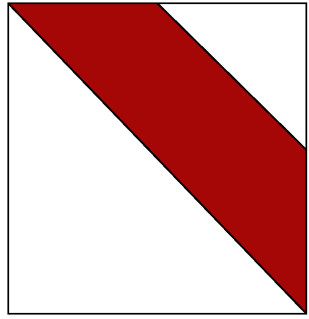}
\caption{A band above the diagonal}\label{fig:band}
\end{center}
\end{figure}

Let $D\in L^\infty([0,1])$ be the identity map from $[0,1]$ to itself;
thus, $D$ generates $L^\infty([0,1])$.
In this paper, with $H$ as above, we compute the free entropy dimension
$\delta_0(Z_H:D)$ of $Z_H$ in the presence of $D$,
in the case $H$ satisfies certain additional hypothesis,
showing that then
\begin{equation}\label{eq:delta0ZH}
\delta_0(Z_H:D)=1+2\,\mathrm{area}(\mathrm{supp}(H)),
\end{equation}
where $\mathrm{supp}(H)$ is the measurable support of $H$
and where the area is Lebesgue measure.
We prove the upper bound $\le$ in~\eqref{eq:delta0ZH} for general $H$,
(see Theorem~\ref{main2}) using basic estimates inspired by~\cite{vN42}.
We prove the lower bound $\ge$ in~\eqref{eq:delta0ZH} for all $H$ 
that are supported in the upper triangle as drawn in Figure~\ref{figure1}
and whose restrictions to some band as drawn in Figure~\ref{fig:band}
are nonzero constant.
(Actually, somewhat weaker conditions suffice --- see Theorem~\ref{main}.)
Our proof of the lower bound uses 
techniques similar to those used in~\cite{DT}.

The organization of the rest of this paper is as follows.
In~\S\ref{sec:DefPrelims}, we discuss some definitions and results
that we need for the calculation.  These include
(\S\ref{subsec:Linf})
basic facts about the class of $L^\infty([0,1])$--circular operators that we consider,
their construction in $L(\F_2)$ and a lemma about them;
(\S\ref{subsec:microT})
a result about certain matrix approximants to the quasinilpotent DT--operator
which was lifted from~\cite{DT} but that follows directly
from work of Aagaard and Haagerup~\cite{Aa} and \'Sniady~\cite{sn};
(\S\ref{subsec:packing}) Jung's equivalent approach to free entropy dimension
in terms of packing numbers~\cite{free};
(\S\ref{subsec:Dyson}) Dyson's formula for the volumes of sets of matrices that
are invariant under unitary conjugation.
In~\S\ref{sec:fedbounds}, we prove the main result, namely
the equation~\eqref{eq:delta0ZH}.
Finally, in~\S\ref{sec:concl}, we consider an example when $\delta_0(Z_H:D)<\delta_0(Z_H)$
and we ask a natural question.
\smallskip
\noindent
{\it Acknowledgement:}  The first named author thanks Kenley Jung for helpful comments.

\section{Definitions and preliminaries}
\label{sec:DefPrelims}

\subsection{$L^{\infty}([0,1])$-circular operators in free group factors}
\label{subsec:Linf}

In this section we recall how 
$L^{\infty}([0,1])$-circular operators in a certain class
were constructed in \cite{B-circ}, and we prove a lemma.
We work in  W$^*$--noncommutative probability space $(\Mcal,\tau)$, with $\tau$ a faithful trace,
and we fix a copy $\Ac=L^\infty[0,1]\subseteq\Mcal$,
such that the restriction of $\tau$ to $\Ac$ is given by integration with respect to
Lebesgue measure on $[0,1]$.
Let $D\in\Ac$ be the operator corresponding the function in $L^\infty[0,1]$ that is the identity
map from $[0,1]$ to itself.
Let $E:\Mcal\to\Ac$ be the $\tau$--preserving conditional expectation.
Let $H\in L^{1}([0,1]^2)$, $H\geq 0$, and assume $H$ has essentially bounded 
coordinate expectations $CE_1(H)$ and $CE_2(H)$, given by 
\begin{equation}\label{eq:CE}
CE_1(H)(x)=\int_0^1{H(x,y)dy},\qquad
CE_2(H)(y)=\int_0^1{H(x,y)dx}.
\end{equation}
By $Z_H$, we will denote an $\Ac$--circular operator in $(\Mcal,E)$ with covariance
$(\alpha_{H},\beta_{H})$
where $\alpha_{H},\,\beta_{H}:L^{\infty}[0,1]\to L^{\infty}[0,1]$ are
given by
\begin{equation}\label{eq:alphaHbetaH}
\alpha_{H}(f)(x)=\int_{0}^{1}{H(t,x)f(t)dt},\qquad\beta_{H}(f)(x)=\int_{0}^{1}{H(x,t)f(t)dt}.
\end{equation}

Suppose $Z\in\Mcal$ is a $(0,1)$--circular element, namely
a circular element  satisfying $\tau(Z)=0$ and $\tau(Z^*Z)=1$,
and suppose $\Ac$ and $\{Z\}$ are $*$--free.
We will construct our operator $Z_H$ from $\Ac$ and $Z$
as in Theorem 6.5 of~\cite{B-circ}.
(Note that our notation differs slightly from that used in~\cite{B-circ}.)

\begin{df} \label{defX}
Let $\omega\in L^{\infty}([0,1]^2)$. We say that $\omega$ is in regular
block form if $\omega$ is constant on all blocks in the regular $n\times n$
lattice
superimposed on $[0,1]^2$, for some $n$, i.e. if there are $n\in\N$ and
$\omega_{i,j}\in\C$, ($1\leq i,j\leq n$) such that
$\omega(s,t)=\omega_{i,j}$ whenever
$\frac{i-1}{n}\leq s \leq \frac{i}{n}$ and $\frac{j-1}{n}\leq t \leq
\frac{j}{n}$, for all integers $1\leq i,j\leq n$.
(We then say $\omega$ is in $n\times n$ regular block form.)
Then we set
\[
M(\omega,Z)=\sum_{i,j=1}^{n}{\omega_{i,j}p_iZp_j}
\]
where $p_i=\mathbf{1}_{[\frac{i-1}{n},\frac{i}{n}]}\in\Ac$.
Note that we have $M(\omega,Z)\in W^{*}(\Ac\cup\{Z\})\cong L(\F_3)$.
\end{df}

Recalling Lemma 6.4 and Theorem 6.5 of~\cite{B-circ} we can state the
following theorem.
\begin{teo}\label{H}
Let $\omega=\sqrt{H}$. Then there
exists a sequence
$\{\omega^{(n)}\}_n$ in $L^{\infty}([0,1]^2)$ such that
\begin{enumerate}
\item[(i)] for each $n$, $\omega^{(n)}$ is in regular block form,
\item[(ii)] $\lim_{n}{\norm{\omega-\omega^{(n)}}_{L^2}}=0$
\item[(iii)] letting $H^{(n)}=(\omega^{(n)})^2$, both
$\norm{CE_1(H^{(n)})}_{\infty}$ and $\norm{CE_2(H^{(n)})}_{\infty}$ remain
bounded as $n$ goes to $\infty.$
\end{enumerate}
Moreover, there is an 
an $L^{\infty}[0,1]$--circular operator $Z_H$ with 
covariance $(\alpha_{H},\beta_{H})$ as described in equations~\eqref{eq:alphaHbetaH}
such that whenever $\{\omega^{(n)}\}_n$ is a sequence satisfying conditions (i)--(iii) above,
the operators $M(\omega^{(n)},Z)$ as given in Definition \ref{defX} converge in
the strong--operator--topology as $n\to\infty$ to $Z_H$.
\end{teo}

\begin{obs}\label{rem:qnDT}
Of particular interest is the operator $Z_R$ when
$R=1_{\{(s,t)\mid s<t\}}$ is the characteristic function in the upper
triangle in $[0,1]^2$.
This $Z_R$ is an instance of the
$\mathrm{DT}(\delta_0,1)$-operator, also called the quasinilpotent DT--operator, and also denoted $T$.
The construction of $Z_R$ in Theorem~\ref{H} above
is approximately what was done in~\S4 of~\cite{DT0}.
\end{obs}

The following lemma will be used in~\S\ref{sec:fedbounds} to prove the upper bound on free
entropy dimension.
For emphasis, we will denote by $\lambda:L^\infty[0,1]\to\Mcal$ the identification of $L^\infty[0,1]$
(with its trace given by Lebesgue measure) and $\Ac=\lambda(L^\infty[0,1])\subseteq\Mcal$.

\begin{lem}\label{lem:Zut}
Let $T=Z_R\in W^*(\{Z\}\cup\Ac)$
be the quasinilpotent DT--operator as described in Remark~\ref{rem:qnDT}.
Let $N$ be an integer, $N\ge2$.
Assume for all $i,j\in\{1,\ldots,N\}$ with $i\ne j$, $Y_{i,j}\in\Mcal$ is a $(0,1)$--circular
element such that the family
\[
\Ac,\quad\{Z\},\quad(\{Y_{i,j}\})_{1\le i,j\le N,\,i\ne j}
\]
is $*$--free.
Let $(e_{ij})_{1\le i,j\le N}$ be a system of matrix units for $M_N(\Cpx)$.
Consider the $*$--noncommutative probability space $(\Mcal\otimes M_N(\Cpx),\tau\otimes\tr_N)$,
and let $\lambdat:L^\infty[0,1]\to\Mcal\otimes M_N(\Cpx)$
be the $*$--homomorphism given by
\[
\lambdat(f)=\sum_{j=1}^N\lambda(f\circ\rho_j)\otimes e_{jj},
\]
where $\rho_j:[0,1]\to[0,1]$ is $\rho_j(t)=\frac tN+\frac{j-1}N$.
Let $\Act=\lambdat(L^\infty[0,1])$.
Then the $\tau\otimes\tr_N$--preserving conditional expectation
$\Et:\Mcal\otimes M_N(\Cpx)\to\Act$ is given by
\[
\Et(\sum_{1\le i,j\le N}a_{ij}\otimes e_{ij})=\sum_{j=1}^NE(a_{jj})\otimes e_{jj}.
\]
Let $c_{ij}\in[0,\infty)$ ($1\le i,j\le N$, $i\ne j$) and let
\[
\Yt=\frac1{\sqrt N}\biggl(\sum_{k=1}^N T\otimes e_{kk}
+\sum_{1\le i,j\le N,\,i\ne j}c_{ij}Y_{ij}\otimes e_{ij}\biggr).
\]
Then $\Yt$ is $\Act$--circular with covariance $(\alpha_H,\beta_H)$ as given
in~\eqref{eq:alphaHbetaH}, where
\[
H(s,t)=
\begin{cases}
1,&\frac{k-1}N\le s\le t\le\frac kN,\,1\le k\le N \\[0.5ex]
(c_{ij})^2,&\frac{i-1}N\le s\le\frac iN,\,\frac{j-1}N\le t\le\frac jN,\,1\le i,j\le N,\,i\ne j.
\end{cases}
\]
\end{lem}
\begin{proof}
Let
\[
\Zt=\frac1{\sqrt N}\biggl(\sum_{k=1}^NZ\otimes e_{kk}
+\sum_{1\le i,j\le N,\,i\ne j}Y_{ij}\otimes e_{ij}\biggr).
\]
We will show that $\Zt$ is $(0,1)$--circular and is $*$--free from $\Act$.
Let $u_1,\ldots,u_N\in\Mcal$ be Haar unitary elements such that the family
\[
(\{u_k,u_k^*\})_{1\le k\le N},\quad\Ac,\quad\{Z\},\quad(\{Y_{i,j}\})_{1\le i,j\le N,\,i\ne j}
\]
is $*$--free (after enlarging $(\Mcal,\tau)$ if necessary).
Let
\[
U=\sum_{k=1}^Nu_k\otimes e_{kk}.
\]
It will suffice to show that $U^*\Zt U$ is $(0,1)$--circular and is $*$--free from $U^*\Act U$.
For this, by results following directly from
Voiculescu's matrix model~\cite{V91} (see~\cite{V90}), it will suffice to show
that each $u_k^*Zu_k$ and each $u_i^*Y_{ij}u_j$ is circular and that the family
\begin{equation}\label{eq:UZU}
(\{u_k^*Zu_k\})_{1\le k\le N},\quad(\{u_i^*Y_{ij}u_j\})_{1\le i,j\le N,\,i\ne j},\quad
(u_k^*\Ac u_k)_{1\le k\le N}
\end{equation}
is $*$--free in $(\Mcal,\tau)$.
Let $Z=V|Z|$ and $Y_{ij}=V_{ij}|Y_{ij}|$ be the polar decompositions.
Then (see~\cite{V90}), $V$ and $V_{ij}$ are Haar unitaries, $|Z|$ and $|Y_{ij}|$ are
quarter--circular elements, $V$ and $|Z|$ are $*$--free
and, for each $i$ and $j$, $V_{ij}$ and $|Y_{ij}|$ are $*$--free in $(\Mcal,\tau)$.
We have the polar decompositions
\begin{align*}
u_k^*Zu_k&=(u_k^*Vu_k)(u_k^*|Z|u_k) \\
u_i^*Y_{ij}u_j&=(u_i^*V_{ij}u_j)(u_j^*|Y_{ij}|u_j).
\end{align*}
Therefore, in order to show that $*$--freeness of the family~\eqref{eq:UZU} and circularity
of $u_k^*Zu_k$ and $u_i^*Y_{ij}u_j$, it will suffice to show $*$--freeness of the family
\[
\begin{gathered}
(\{u_k^*|Z|u_k\})_{1\le k\le N},\quad(\{u_k^*Vu_k\})_{1\le k\le N}, \\
\quad(\{u_j^*|Y_{ij}|u_j\})_{1\le i,j\le N,\,i\ne j},\quad(\{u_i^*V_{ij}u_j\})_{1\le i,j\le N,\,i\ne j},\quad
(u_k^*\Ac u_k)_{1\le k\le N}.
\end{gathered}
\] 
Let $B$ be a Haar unitary generating $W^*(|Z|)$, let $B_{ij}$ be a Haar unitary generating $W^*(|Y_{ij}|)$,
and let $C$ be a Haar unitary generating $\Ac$.
It will suffice to show $*$--freeness of the family
\[ 
\begin{gathered}
(u_k^*Bu_k)_{1\le k\le N},\quad(u_k^*Vu_k)_{1\le k\le N}, \\
\quad(u_j^*B_{ij}u_j)_{1\le i,j\le N,\,i\ne j},\quad(u_i^*V_{ij}u_j)_{1\le i,j\le N,\,i\ne j},\quad
(u_k^*Cu_k)_{1\le k\le N}
\end{gathered}
\] 
of Haar unitaries.
This follows from the $*$--freeness of the family
\[
B,\,C,\,V,\,(u_k)_{1\le k\le N},\quad(B_{ij})_{1\le i,j\le N,\,i\ne j},\quad
(V_{ij})_{1\le i,j\le N,\,i\ne j}.
\]
by an argument involving words in a free group.
This shows that $\Zt$ is $(0,1)$--circular and $*$--free from $\Act$.

Now we use the method of Theorem~6.5 of~\cite{B-circ}, described in Theorem~\ref{H}
above, but taking $\omega^{(n)}$ in $n\times n$ regular block form with $n$ always
a multiple of $N$, and with each such $\omega^{(n)}$ constant equal to $c_{ij}$ 
on each off--diagonal block of the form
$[\frac{i-1}N,\frac iN]\times[\frac{j-1}N,\frac jN]$ for $1\le i,j\le N$,
$i\ne j$,
where projections from $\Act$ are used to cut $\Zt$ and make each $M(\omega^{(n)},\Zt)$.
It is then clear that the operators $M(\omega^{(n)},\Zt)$ converge to $\Yt$ as $n\to\infty$,
and, from Theorem~\ref{H}, they also converge to an $\Act$--circular operator having the desired
covariance $(\alpha_H,\beta_H)$.
\end{proof}

\subsection{Microstates for the quasinilpotent DT-operator}
\label{subsec:microT}

Let $T=Z_R$ be the quasinilpotent DT--operator as described in Remark~\ref{rem:qnDT}
and let $D$ be the corresponding operator described in~\S\ref{subsec:Linf}.
It was proved by Aagard and Haagerup \cite{Aa} that if we consider
$T$ a DT$(\delta_0,1)$-operator and $Y$ a circular operator that is $*$--free from $T$ (and $D$),
then the Brown measure of $T+\epsilon Y$ is equal to the uniform distribution
on the closed disk centered at 0 and of radius
$r_{\epsilon}=\log(1+\epsilon^{-2})^{-\frac{1}{2}}$. Note how
slowly this disk shrinks as $\epsilon$ approaches to 0. Moreover,
they also showed that the spectrum of $T+\epsilon Y$ is equal to the
disk.

The next lemma is an immediate consequence of the above described Brown measure result
of Aagaard and Haagerup and
a result of \'Sniady \cite{sn}.
A detailed proof can be formulated exactly as was done for Lemma 2.2 in \cite{DT}.
In
the following lemma and throughout this paper, for a matrix $A\in
M_{k}(\C)$ we let $|A|_2=\tr_{k}(A^*A)^{1/2}$, where
$\tr_k$ is the normalized trace of $M_{k}(\C)$. Also, by the
eigenvalue distribution of a matrix $A\in M_{k}(\C)$ we mean
the probability measure $\frac1n\sum_1^n\delta_{\lambda_j}$,
where
$\lambda_1,\ldots,\lambda_k$ are the eigenvalues of $A$ listed according to
general multiplicity.

\begin{lem}\label{lemma}
Let $c>0$. Then there exists sequences $\{g_k\}_k$ and $\{y_k\}_k$
such that for
any $\epsilon>0$, there exists a sequence $\{z_{k,\epsilon}\}_k$
such that
\begin{itemize}
\item $g_k,y_k,z_{k,\epsilon}\in M_{k}(\C)$,
\item $\norm{g_k}$, $\norm{y_k}$ and $\norm{z_{k,\epsilon}}$ remain bounded as
$k\to +\infty$,
\item $\limsup_{k}{|y_k-z_{k,\epsilon}|_2}\leq \epsilon c$,
\item the pair $(g_k,y_{k})$ converges in $*$--moments as $k\to +\infty$ to the pair $(D,T)$,
\item the eigenvalue distribution of $z_{k,\epsilon}$ converges
weakly as $k\to +\infty$ to the measure $\sigma_{\epsilon,c}$,
which is the uniformly distributed measure in
the disk of center at 0 and radius
$r_{\epsilon,c}=c\log(1+\epsilon^{-2})^{-\frac{1}{2}}$ in the
complex plane.
\end{itemize}
\end{lem}

\subsection{Packing number formulation of the free entropy dimension}
\label{subsec:packing}

In this section we will review the packing number formulation of
Voiculecu's microstates free entropy dimension due to K. Jung
\cite{free}. If $X=(x_1,\ldots,x_{n})$ and $Z=(z_1,\ldots,z_m)$ are
tuples of
selfadjoint elements in a tracial von Neumann algebra, then the
microstates free entropy dimension (as defined by Voiculescu
\cite{V}) is given by the formula
$$\delta_{0}(X)=n+\limsup_{\epsilon\to 0}{\frac{\chi(x_1+\epsilon
s_1,\ldots,x_n+\epsilon s_n:s_1,\ldots,s_n)}
{|\log \epsilon|}}$$
and the microstates free entropy dimension in the presence of $Z$ is defined
by
$$\delta_{0}(X:Z)=n+\limsup_{\epsilon\to 0}{\frac{\chi(x_1+\epsilon
s_1,\ldots,x_n+\epsilon s_n:z_1,\ldots,z_m,s_1,\ldots,s_n)}
{|\log \epsilon|}}$$
where $\{s_1,\ldots,s_n \}$ is a semicircular
family free from $X$ and $Z$. The packing formulation found in \cite{free}
is
\begin{equation}\label{pack}
\delta_{0}(X)=\limsup_{\epsilon\to
0}{\frac{\mathbb{P}_{\epsilon}(X)}{|\log
\epsilon|}}\,\,\,\,\,\,\,\,\,\,\,\,\,\,\,\,\,\,\,\,\,\delta_{0}(X:Z)=\limsup_{\epsilon\to
0}{\frac{\mathbb{P}_{\epsilon}(X:Z)}{|\log \epsilon|}}
\end{equation}
where
$$\mathbb{P}_{\epsilon}(X)=\inf_{m,\gamma}\limsup_{k}k^{-2}\log
P_{\epsilon}(\Gamma(X;m,k,\gamma))$$
and
$$\mathbb{P}_{\epsilon}(X:Z)=\inf_{m,\gamma}\limsup_{k}k^{-2}\log
P_{\epsilon}(\Gamma(X:Z;m,k,\gamma))$$
Here, $\Gamma(X:Z;m,k,\gamma)\subseteq (M_{k}(\C)_{s.a.})^{n}$ is
the microstates space of Voiculescu, and $P_{\epsilon}$ is the
packing number with respect to the metric arising from the
normalized trace.
Let $Y=(y_1,\ldots,y_n)$ and $W=(w_1,\ldots,w_m)$
be arbitrary
tuples of possibly non-selfadjoints elements in a tracial von
Neumann algebra. Now the definition of $\mathbb{P}_\epsilon$ makes
perfect sense for the set $Y$ if we replace the microstates space
in (\ref{pack}) with the non-selfadjoint $*-$microstates space
$\Gamma(Y:W;m,k,\gamma)\subseteq (M_{k}(\C))^{n}$, which is the set
of all $n-$tuples of $k\times k$ matrices whose $*-$moments up to
order $m$ approximate those of $Y$ within tolerance of $\gamma$ in the
presence of $W$.
It is also true that
$$\delta_{0}(\text{Re}(y_1),\text{Im}(y_1),\ldots,\text{Re}(y_n),\text{Im}(y_n):W)=\limsup_{\epsilon\to
0}{\frac{\mathbb{P}_{\epsilon}(Y:W)}{|\log \epsilon|}}$$
see \cite{DT} for details.

Finally, we review the standard volume comparison inequality for packing numbers.
Recall that for a metric space $A$ we have
\[
P_{4\eps}(A)\le K_{2\eps}(A)\le P_\eps(A),
\]
where $P_\eps(A)$ is the $\eps$--packing number, i.e.\ the maximal number of disjoint open
balls of radius $\eps$ in $A$, and $K_\eps(A)$ is the minimal number of elements in a cover of $A$
consisting of open balls of radius $\eps$.
If $A$ is a subspace of a Euclidean space, then we have 
\[
\vol(\Nc_\eps(A))\le K_\eps(A)\cdot\vol(\Bc_{2\eps}),
\]
where $\Nc_\eps(A)$ is the $\eps$--neighborhood,  $\Bc_r$ is
a ball of radius $r$ and $\vol$ is the volume, all in the ambient Euclidean space.
We thus have the volume comparison test,
\begin{equation}\label{eq:vct}
P_\eps(A)\ge K_{2\eps}(A)\ge\frac{\vol(\Nc_{2\eps}(A))}{\vol(\Bc_{4\eps})}.
\end{equation}

\subsection{Dyson's formula}
\label{subsec:Dyson}

Every matrix of $M_{k}(\C)$ has an upper-triangular matrix in its
unitary orbit. Thus, letting $T_{k}(\C)$ denote the set of
upper-triangular matrices in $M_{k}(\C)$, there is a probability measure
$\nu_k$ on $T_{k}(\C)$ such that
$$\lambda_{k}(\mathcal{O})=\nu_{k}(\mathcal{O}\cap T_{k})$$
for every $\mathcal{O}\subseteq M_{k}(\C)$ that is invariant under unitary 
conjugation.
Freeman Dyson identified such a measure \cite{mehta}, and showed
that if we view $T_{k}(\C)$ as a Euclidean space of real dimension
$k(k+1)$ with coordinates corresponding to the real and imaginary
part of the matrix entries lying on and above the diagonal, then
$\nu_{k}$ is absolutely continuous with respect to Lebesgue
measure on $T_{k}(\C)$ and has density given at $A=(a_{ij})_{1\leq
i,j \leq k}\in T_{k}(\C)$ by
\begin{equation}
C_{k}\cdot\prod_{1\leq p < q\leq k}{|a_{pp}-a_{qq}|^2}
\,\,\,\,\,\text{where}\,\,\,\,\,C_{k}=\frac{\pi^{k(k+1)/2}}{\prod_{j=1}^{k}{j!}}.
\end{equation}
We will use Dyson's formula in our main result to find lower bound
on the volume of unitary orbits of an $\epsilon-$neighborhood of
the microstates space.


\section{Free entropy dimension computations}
\label{sec:fedbounds}

\begin{lem}\label{prop1}
Let $(\Omega,\mu)$ a finite measurable space. Let $f\in
L^1(\Omega)$ and $f\geq 0$.
Then
$$\lim_{\epsilon\to
0}{\frac{\int_{\Omega}{\log(\max(f(t),\epsilon))}d\mu(t)}{|\log \epsilon|}}=
\mu(\mathrm{supp}(f))-\mu(\Omega),$$
where $\mathrm{supp}(f)=f^{-1}((0,+\infty))$.
\end{lem}
\begin{proof}
It is clear that we have
$\log(\max(f(t),\epsilon))\leq
\log(f(t)+1)+\log(\epsilon)\cdot\mathbf{1}_{f^{-1}([0,\epsilon))}$,
and this yields
$$\limsup_{\epsilon\to
0}{\frac{\int_{\Omega}{\log(\max(f(t),\epsilon))}d\mu(t)}{|\log \epsilon|}}
\leq -\liminf_{\epsilon\to
0}{\mu(f^{-1}([0,\epsilon)))}=\mu(\mathrm{supp}(f))-\mu(\Omega).$$
On the other hand, given $\gamma>0$, let $\delta>0$ be such that
$\mu(f^{-1}((0,\delta)))<\gamma$.
Taking $0<\epsilon<\delta$,
we have $\mathbf{1}_{f^{-1}([0,\delta))}\cdot
\log\epsilon+\mathbf{1}_{f^{-1}([\delta,+\infty))}\cdot\log\delta
\leq \log{\max(f(t),\epsilon)}$ and integrating on both sides we
obtain
$$\mu(f^{-1}([0,\delta)))\cdot\log\epsilon+\mu(f^{-1}([\delta,+\infty)))\cdot\log\delta\leq
\int_{\Omega}{\log(\max(f(t),\epsilon))}d\mu(t).$$ Now dividing by
$|\log\epsilon|$ and taking $\liminf$ on both sides we get
$$-\mu(f^{-1}([0,\delta)))\leq \liminf_{\epsilon\to
0}{\frac{\int_{\Omega}{\log(\max(f(t),\epsilon))}d\mu(t)}
{|\log \epsilon|}}.$$ Using the fact that
$\mu(f^{-1}([0,\delta)))<\mu(f^{-1}(0))+\gamma$ and
that $\gamma$ is arbitrary we obtain 
$$\mu(\mathrm{supp}(f))-\mu(\Omega)\leq \liminf_{\epsilon\to
0}{\frac{\int_{\Omega}{\log(\max(f(t),\epsilon))}d\mu(t)}{|\log
\epsilon|}},$$
proving the claim.
\end{proof}

As in~\S\ref{subsec:Linf}, we work in $(\Mcal,\tau)$ and
we have $\Ac=L^{\infty}[0,1]$ and a $(0,1)-$circular element $Z$
such that $\Ac$ and $Z$ are $*-$free, and with $H$ as described there.
We construct as in~\S\ref{subsec:Linf}
an $L^{\infty}[0,1]$--circular operator $Z_H\in W^*(\Ac\cup\{Z\})\cong L(\F_3)$.
We also take $D=D^{*}\in\Ac$ to correspond to the identity function from $[0,1]$ to itself.
The following is our main result.
\begin{teo}\label{main}
Let $H\geq 0$, $H\in L^{1}([0,1]^2)$ have essentially bounded coordinate expectations $CE_1(H)$
and $CE_2(H)$, as in equations~\eqref{eq:CE}.
Assume $H$ has support contained in the upper-triangle $U$ of $[0,1]^2$ and assume there exists
$r\in\N$ such that 
\[
\Delta:=\bigcup_{i=1}^{r}U_i^{(r)}\subseteq
\mathrm{supp(H)},\qquad U_{i}^{(r)}=\{(x,y):\frac{i-1}{r}\leq 
x<y\leq \frac{i}{r}\}
\]
and that $H$ restricted to $\Delta$ is constant equal
to $c>0$. Then
\[
\delta_0(Z_{H}:D)\geq 1+2\cdot\mathrm{area}(\mathrm{supp}(H)).
\]
In particular,
$\delta_0(Z_{H})\geq 1+2\cdot\mathrm{area}(\mathrm{supp}(H))$.
\end{teo}

\begin{proof}
Without loss of generality we can assume $c=1$. Fix $\epsilon>0$.
By hypothesis we may choose $N$ arbitrarily large and so that
$\bigcup_{i=1}^{N}U_{i}^{(N)}\subseteq \Delta$.
Let $R>1$, $m\in\N$ and $\gamma>0$.
There is $\delta>0$ such that $\|Z_H-Y\|_2<\delta$ implies
$\Gamma_{R}(Y;m,k,\gamma/2)\subseteq\Gamma_{R}(Z_{H};m,k,\gamma)$.
Making use of Theorem~\ref{H},
there exist $M=Np$ and
$$\omega:=\sum_{i=1}^{M}{\mathbf{1}_{U_i^{(M)}}}+\sum_{1\leq i<j\leq
M}{\alpha_{ij}\mathbf{1}_{E_{ij}^{(M)}}}$$ where
$E_{ij}^{(M)}=\{(x,y):\frac{i-1}{M}\leq x\leq
\frac{i}{M}\,,\,\frac{j-1}{M}\leq y\leq \frac{j}{M}\}$ with
$\alpha_{ij}>0$, such that
$\norm{Z_{H}-Z_{\omega}}_2<\delta$ and, therefore, we have
$\Gamma_{R}(Z_{\omega};m,k,\gamma/2)\subseteq
\Gamma_{R}(Z_{H};m,k,\gamma)$.
We define the sets of indices $$\Theta=\{(i,j):1\leq
i<j\leq p,\,\,p+1\leq i<j\leq 2p,\ldots,(N-1)p+1\leq i<j\leq
Np\}$$ and
$$\Phi=\{(i,j):1\leq i<j\leq Np\}\setminus\Theta.$$
For example, in the case $N=2$ and $p=4$
the squares corresponding to elements of $\Theta$ are shaded
in Figure~\ref{fig:thetashaded}.
\begin{figure}[Ht]\label{fig:thetashaded}
\begin{center}
\includegraphics[width=6cm]{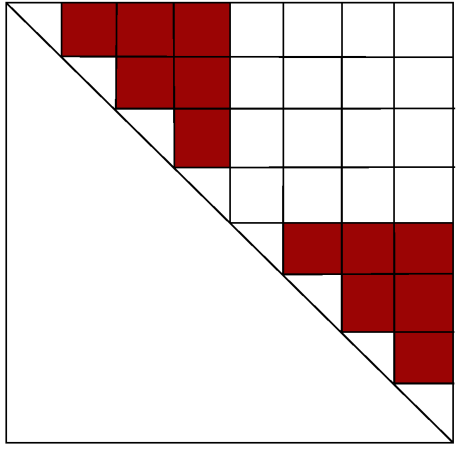}
\caption{Case $N=2$ and $p=4$}
\end{center}
\end{figure}
Note that by the hypothesis of $H$ we may insist, $\alpha_{ij}=1$
whenever $(i,j)\in\Theta$. Let $\gamma'=\gamma/(MR)^{m-1}$.

Consider $(C_{11},\ldots,C_{MM}),\,(C_{ij})_{1\leq i<j\leq M}$ a
$*-$free
family in $(\Mcal,\tau)$, where
each $C_{ii}$ is $\DT(\delta_0,\frac{1}{\sqrt{M}})$, and each $C_{ij}$ with
$i<j$
is circular with $\tau(|C_{ij}^2|)=\frac{1}{M}$.
Let $\{g_k\}_k$ and $\{y_k\}_k$ the sequences constructed in Lemma~\ref{lemma} with $c=1/\sqrt M$.
There are
$a_{ij}(k)\in M_{k}(\C)$ for $(i,j)\in\Theta$ such that
for each $(i,j)\in\Theta$ as before $a_{ij}(k)$ converge in distribution as
$k\to +\infty$ to a $(0,\frac{1}{M})$-circular element and such
that the family
\[
\{g(k),y(k)\},\,(\{a_{ij}(k)\})_{(i,j)\in\Theta}
\]
of sets of random variables is asymptotically $*$--free as $k\to\infty$.
By an application
of Corollary 2.14 of \cite{Voi}, for $k$ large enough there exists
a set $\Omega_k\subset
\Gamma((C_{ij})_{(i,j)\in\Phi};m,k,\gamma')$ such that for any
$(\eta_{ij})_{(i,j)\in\Phi}\in \Omega_k$,
\[
\{y_k,g(k)\},\,(a_{ij}(k))_{(i,j)\in\Theta},\,(\eta_{ij})_{(i,j)\in\Phi}
\]
is an $(m,\gamma')$--$*$--free family of sets of random variables and
\begin{multline}\label{eqx}
\liminf_{k}{\Bigg(k^{-2}.\log(\text{vol}(\Omega_k))+\bigg(\frac{N(N-1)p^2}{2}\bigg).\log(k)\Bigg)}
\geq \\
\geq \chi((\RePt C_{ij})_{(i,j)\in\Phi},(\ImPt C_{ij})_{(i,j)\in
\Phi})>-\infty
\end{multline}
where the volume is computed with respect to the Euclidean norm
$k^{1/2}|\cdot|_2$. For each  $(\eta_{ij})_{(i,j)\in\Phi}\in
\Omega_k$ we define a matrix $R(k)\in M_{Mk}(\C)$ by
\[ 
R(k)=\left[ \begin{array}{cccc}
r_{11}(k)& r_{12}(k) & \ldots  & r_{1M}(k) \\
0 & r_{22}(k) & \ldots & r_{2M}(k)\\
\vdots & \ddots & \ldots & \vdots\\
0 & \ldots & 0 & r_{MM}(k)\\
\end{array}\right],
\qquad
r_{ij}(k)=\begin{cases}
y_k,&i=j \\
a_{ij},&(i,j)\in \Theta \\
\alpha_{ij}\eta_{ij},&(i,j)\in\Phi.
\end{cases}
\]
Let
\[
G(k)=\text{diag}(g(k),{\textstyle\frac{1}{M}}+g(k),\ldots,{\textstyle\frac{M-1}{M}}+g(k))\in M_{Mk}(\C).
\]
As a consequence of Lemma~\ref{lem:Zut},
\[
(R(k),G(k))\in \Gamma(Z_{\omega},D;m,Mk,\gamma/2).
\]
Set $\tilde{\alpha}_{ij}=\max(\alpha_{ij},\epsilon)$ and let
\begin{equation}\label{eq:Rtil}
\tilde{R}(k)=\left[ \begin{array}{cccc}
r_{11}(k)& r_{12}(k) & \ldots  & r_{1M}(k) \\
0 & r_{22}(k) & \ldots & r_{2M}(k)\\
\vdots & \ddots & \ldots & \vdots\\
0 & \ldots & 0 & r_{MM}(k)\\
\end{array}\right],
\qquad
r_{ij}(k)=\begin{cases}
y_k,&i=j\\
a_{ij},&(i,j)\in \Theta\\
\tilde{\alpha}_{ij}\eta_{ij},&(i,j)\in \Phi.\\
\end{cases}
\end{equation}
Then $\tilde{R}(k)$ lies in an $\epsilon$-neighborhood of
$\Gamma(Z_{\omega}:D;m,Mk,\gamma/2)$.
Let $A_{l}(k)\in M_{kp}(\C)$ for $l\in \{1,2,\ldots,N\}$ be defined by
\[
A_{l}(k)=\left[ \begin{array}{cccc}
y_k & a_{f+1,f+2} & \ldots  & a_{f+1,f+p}\\
0 & y_k & \ldots &  \vdots\\
\vdots & \ddots & \ldots & a_{f+p-1,f+p}\\
0 & \ldots & 0 & y_k\\
\end{array} \right]
\]
with $f=(l-1)p$.
Note that we have
\begin{equation}\label{eq:RtilD}
\tilde{R}(k)=
\left[ \begin{array}{cccc}
A_{1}(k)& Y_{12}(k)& \ldots  & Y_{1N}(k)\\
0 & A_{2}(k)& \ldots &\vdots\\
\vdots& \ldots & \ddots & Y_{N-1,N}\\
0 & \ldots & 0 & A_{N}(k)\\
\end{array} \right],
\end{equation}
where the $Y_{ij}(k)\in M_{pk}(\C)$ are determined by equations~\eqref{eq:Rtil} and~\eqref{eq:RtilD}.
Then, by again making use of Lemma~\ref{lem:Zut}, we have $A_{l}(k)\in
\Gamma_{p^2R}(\frac1{\sqrt N}T;m,pk,\gamma)$ for all $l\in
\{1,2,\ldots,N\}$, where $T$ is the the $\DT(\delta_0,1)$--operator.
Let $\epsilon>0$ and let $z_{k,\epsilon}$ be as in Lemma~\ref{lemma}.
Let
\begin{displaymath}
B_{l,\epsilon}(k)=\left[ \begin{array}{cccc}
z_{k,\epsilon} & a_{f+1,f+2} & \ldots  & a_{f+1,f+p}\\
0 & z_{k,\epsilon} & \ldots &  \vdots\\
\vdots & \ddots & \ldots & a_{f+p-1,f+p}\\
0 & \ldots & 0 & z_{k,\epsilon}\\
\end{array} \right]\in M_{kp}(\C).
\end{displaymath}
Note that the eigenvalue distribution of $B_{l,\epsilon}(k)$
converge weakly as $k\to
+\infty$ to the measure $\sigma_{\epsilon,\frac1{\sqrt{N}}}$ of Lemma \ref{lemma}.

Since every complex matrix can be put into an upper-triangular form with
respect to an orthonormal basis, we can find a
$k\times k$ unitary
matrix $v(k)$ such that $v(k)z_{k,\epsilon}v(k)^{*}$ is upper
triangular.
Since microstate spaces are invariant under conjugation by unitaries,
also $(v(k)\otimes I_M)\tilde{R}(k)(v(k)\otimes I_M)^*$
lies in an $\epsilon$-neighborhood of
$\Gamma(Z_{\omega}:D;m,Mk,\gamma/2)$.

For each $1\leq l\leq N$, we have
\[
|(v(k)\otimes I_p)B_{l,\epsilon}(k)(v(k)\otimes I_p)^*-(v(k)\otimes I_p)A_{l}(k)(v(k)\otimes I_p)^*|_2
=|A_{l}(k)-B_{l,\epsilon}(k)|_2.
\]
Since $\limsup_{k}|B_{l,\epsilon}(k)-A_{l}(k)|_2\leq
\frac{\epsilon}{\sqrt{N}}$, and taking $N>4$, for $k$
sufficiently large we have
\[
|(v(k)\otimes I_p)B_{l,\epsilon}(k)(v(k)\otimes I_p)^*-(v(k)\otimes I_p)A_{l}(k)(v(k)\otimes I_p)^*|_2
\leq\epsilon/2.
\]
Set $\Bt_{l}(k)=(v(k)\otimes I_p)B_{l,\epsilon}(k)(v(k)\otimes I_p)^*$ and
$\Yt_{ij}(k)=(v(k)\otimes I_p)Y_{ij}(k)(v(k)\otimes I_p)^*$ and denote by
$\Gc_k$ the set of all $Mk\times Mk$ matrices of the form
\begin{displaymath}
\left[ \begin{array}{cccc}
\Bt_1(k)& \Yt_{12}(k)& \ldots  & \Yt_{1N}(k) \\
0 & \Bt_{2}(k) & \ddots  & \ldots \\
\vdots & \ddots & \ddots & \Yt_{N-1,N}(k) \\
0 & \ldots & 0 & \Bt_{N}(k)\\
\end{array} \right],
\end{displaymath}
over all choices of $(\eta_{ij})_{(i,j)\in\Phi}\in\Omega_k$.
Note that the matrices in $\Gc_{k}$ are upper triangular and
their eigenvalue distributions are exactly the same as
$z_{k,\epsilon}$.
For $k$ sufficiently large, the set $\Gc_{k}$ lies in a
$2\epsilon-$neighborhood of $\Gamma(Z_{\omega}:D;m,Mk,\gamma/2)$ 
and, therefore, in a $2\epsilon-$neighborhood of
$\Gamma(Z_{H}:D;m,Mk,\gamma)$.
Let $\theta(\Gc_k)$ denote the unitary orbit of
$\Gc_k$ in $M_{Mk}(\C)$.
We will now find lower bounds for the
volumes of
$\theta(\Gc_k)$ and thus, via the estimate~\eqref{eq:vct}, lower bounds
for packing number of $\Gamma(Z_{H}:D;m,Mk,\gamma)$.

Denote by $\Hc_k\subset M_{Mk}(\C)$ the set of all matrices of the form
\begin{displaymath}
\left[ \begin{array}{cccc}
0& \Yt_{12}(k)& \ldots  & \Yt_{1N}(k)\\
0 & 0 & \ddots  & \ldots \\
\vdots & \ddots & \ddots & \Yt_{N-1,N}(k)\\
0 & \ldots & 0 & 0\\
\end{array} \right],
\end{displaymath}
over all choices of $(\eta_{ij})_{(i,j)\in\Phi}\in\Omega_k$.
Notice that $\Hc_{k}$ is
isometric to the space of all matrices of the form
$(w_{ij})_{1\leq i,j\leq M}\in M_{Mk}(\C)$
with $w_{ij}\in M_{k}(\C)$ and
\[
w_{ij}=\begin{cases}
0, &(i,j)\notin \Phi \\
\tilde{\alpha_{ij}} \eta_{ij}, & (i,j)\in \Phi.
\end{cases}
\]
It follows that
$\Hc_{k}$ must also have the same volume as the above subspace,
computed in the ambient Hilbert space of block upper-triangular matrices
with the indicated entries set to zero.
Therefore,
$$\mathrm{vol}(\Hc_k)=\mathrm{vol}(\Omega_k)\cdot(M^{1/2})^{k^2 M(M-1)}\cdot
\prod_{(i,j)\in \Phi}{|\tilde{\alpha}_{ij}|^{2k^2}}.$$
Let $T_{n}$ the set of upper triangular matrices in
$M_{n}(\C)$; let $T_{n,<}$ denote the matrices in $T_n$ that have zero
diagonal, i.e. the strictly
upper triangular matrices.
Denote by $\Wc_k$ the set of $T_{Mk,<}$ consisting
of
all matrices $x$ such that $|x|_2<\epsilon$ and $x_{ij}=0$ whenever $1\leq
r<s\leq N$
and $(r-1)pk<i\leq rpk$, $(s-1)pk<j\leq spk$. Thus, $\Wc_{k}$ consists of
$N\times N$
diagonal matrices whose diagonal entries are strictly upper triangular
$pk\times pk$ matrices.
Denote by $\Dc_{k}$ the subset of diagonal matrices $x$ of $M_{Mk}(\C)$ such
that $|x|_2<\epsilon$.
It follows that if $f_{k}$ is the matrix
\begin{displaymath}
f_{k}=\left[ \begin{array}{cccc}
\Bt_{1}(k) & 0 & \ldots & 0\\
0 & \Bt_{2}(k) & \ddots & \vdots\\
\vdots & \ddots & \ddots & 0\\
0 & \ldots & 0 & \Bt_{N}(k)\\
\end{array} \right]
\end{displaymath}
then $f_{k}+\Dc_{k}+\Wc_k+\Hc_{k} \subset \mathcal{N}_{3\epsilon}(\Gc_k)$,
where the $3\epsilon-$neighborhood is taken in the ambient space $T_{Mk}$
with respect
to the metric induced by $|\cdot|_2$.
Now observe that the space of diagonal $Mk\times Mk$ and $T_{Mk,<}$ are
orthogonal subspaces.
Let $\theta_{3\epsilon}(\Gc_{k})$
denote the $3\epsilon$ neighborhood of the unitary orbit of $\theta(\Gc_{k})$
of $\Gc_{k}$.
Let $dX$ denote Lebesgue measure on $T_{Mk}$ corresponding to the Euclidean norm
$(Mk)^{1/2}|\cdot|_2$, which is coordinatized by the complex entries
$X=\{x_{ij}\}_{1\leq i\leq j\leq Mk}$ of the matrix.
Using Dyson's formula we have
\begin{equation}\label{eq}
\begin{aligned}
\mathrm{vol}(\theta_{3\epsilon}(\Gc_{k})) & \geq
C_{Mk}\cdot\int_{f_{k}+\Dc_{k}+\Wc_k+\Hc_{k}}
  {\prod_{1\leq i<j\leq Mk}{|x_{ii}-x_{jj}|^2dX}} \\
& =  C_{Mk}\cdot \mathrm{vol}(\Wc_k+\Hc_k)\cdot
\int_{D(f_{k}+\Dc_{k})}{\prod_{1\leq i<j\leq Mk}{|x_{ii}-x_{jj}|^2dx_{11}\cdots dx_{Mk,Mk}}} \\
& \geq C_{Mk}\cdot \mathrm{vol}(\Wc_k+\Hc_k)\cdot E_{\epsilon}(f_{k})
\end{aligned}
\end{equation}
where the constant $C_{Mk}$ is as in \cite{DT}  and where
$\mathrm{vol}(\theta_{3\epsilon}(\Gc_{k}))$ is computed in
$M_{Mk}(\C)$ and $\mathrm{\Wc_k+\Hc_k}$ is computed in $T_{Mk,<}$,
both being Euclidean volumes corresponding to the norms
$(Mk)^{1/2}|\cdot|_2$,
where the integral over $D(f_k+\Dc_k)$ is over the diagonal parts of these matrices,
and where $E_\epsilon(f_k)$ is the integral defined on p.~252 of~\cite{DT}.
It is clear that
$\theta_{3\epsilon}(\Gc_{k})\subset
\mathcal{N}_{4\epsilon}(\Gamma(Z_{H}:D;m,Mk,\gamma))$, so
(\ref{eq}) gives a lower bound on
$\mathrm{vol}(\mathcal{N}_{4\epsilon}(\Gamma(Z_{H}:D;m,Mk,\gamma)))$.\\
Using (\ref{eq}) and the standard volume comparison test~\eqref{eq:vct}, we have
\begin{align*}
P_{2\epsilon}(\Gamma(Z_{H}:D;m,Mk,\gamma)) & \geq
\frac{\mathrm{vol}(\mathcal{N}_{4\epsilon}(\Gamma(Z_{H};m,Mk,\gamma)))}
{\mathrm{vol}(\mathcal{B}_{8\epsilon})}\nonumber \\
& \geq   C_{Mk}\cdot \mathrm{vol}(\Wc_k+\Hc_k)\cdot
E_{\epsilon}(f_{k})\cdot
\frac{\Gamma((Mk)^2+1)}{\pi^{(Mk)^2}(8(Mk)^{1/2}\epsilon)^{2(Mk)^2}}
\end{align*}
where $\mathcal{B}_{8\epsilon}$ is a ball in $M_{Mk}(\C)$ of radius
$8\epsilon$ with respect to $|\cdot|_2$, and we are taking volumes
corresponding to the Euclidean
norm $(Mk)^{1/2}|\cdot|_2$. Since $\Wc_k$ and $\Hc_k$ are orthogonal, we
have
that $\mathrm{vol}(\Wc_k+\Hc_k)=\mathrm{vol}(\Wc_k)\cdot
\mathrm{vol}(\Hc_k)$,
where each volume is taken in the subspace of appropriate dimension. But
$\Wc_k$
is a ball of radius $(Mk)^{1/2}\epsilon$ in a space of real
dimension $Npk(pk-1)$, so
$$\mathrm{vol}(\Wc_k+\Hc_k)=\frac{\pi^{\frac{Npk(pk-1)}{2}}((Mk)^{1/2}\epsilon)^{Npk(pk-1)}}{\Gamma(\frac{Npk(pk-1)}{2}+1)}\cdot
\mathrm{vol}(\Hc_k)$$
where $\mathrm{vol}(\Hc_k)=\mathrm{vol}(\Omega_k)\cdot(M^{1/2})^{k^2
M(M-1)}\cdot \prod_{(i,j)\in \Phi}{|\tilde{\alpha}_{ij}|^{2k^2}}.$
Using Stirling's formula and $M=Np$, we find
\begin{eqnarray*}
\mathbb{P}_{\epsilon}(Z_{H}:D;m,\gamma) & \geq & \liminf_{k}{(Mk)^{-2}\log
P_{\epsilon}(\Gamma(Z_{H}:D;m,Mk,\gamma))} \\
& \geq & \liminf_{k}{(Mk)^{-2}\log(E_{\epsilon}(f_k))} \\
&& + \quad  \liminf_{k}\Bigg( (Mk)^{-2}\log(C_{Mk}) +
(Mk)^{-2}\log(\mathrm{vol}(\Omega_k))+\nonumber\\
&& + \quad \Big(2-\frac{1}{N}\Big)|\log\epsilon|+\Big(1-\frac{1}{2N}\Big)\log k \\
&& + \quad (\frac{M-1}{2M})\log M
+\frac{2}{M^2}\sum_{(i,j)\in\Phi}{\log|\tilde{\alpha_{ij}}|}\Bigg)+
L_1 \\
& = & \liminf_{k}
(Mk)^{-2}\log(E_{\epsilon}(f_k)) \\
&& + \quad\liminf_{k}\Big((Mk)^{-2}\log
C_{Mk}+\frac{1}{2}\log Mk \Big) \\
&& + \quad \liminf_{k} \Bigg ( (Mk)^{-2}\log(\mathrm{vol}(\Omega_k)) +
\Big (\frac{1}{2}-\frac{1}{2N}\Big ) \log k \Bigg ) \\
&& + \quad \Big ( 2-\frac{1}{N}\Big)|\log \epsilon|
+\frac{2}{M^2}\sum_{(i,j)\in\Phi}{\log|\tilde{\alpha_{ij}}|} +
L_2
\end{eqnarray*}
where $L_1$ and $L_2$ are constants independent of $\epsilon, m$
and $\gamma$.
As $\gamma\to 0$ and $m\to +\infty$, we have convergence
$$\frac{2}{M^2}\sum_{(i,j)\in\Phi}{\log|\tilde{\alpha_{ij}}|}\longrightarrow
2\iint_{K_N}\log(\max(H(s,t),\epsilon))dsdt$$
where 
\[
K_N=\bigcup_{j=1}^{N-1}\bigg\{\frac{j}{N}\leq x\leq \frac{j+1}{N}\leq y \leq 1\bigg\}.
\]
Note that we have
$\mathrm{area}(K_N)=\frac{N(N-1)}{2N^2}$.
Now by~\eqref{eqx}, we have
\begin{multline*}
\liminf_{k}{\Bigg((Mk)^{-2}\log(\text{vol}(\Omega_k))
 +\bigg(\frac{1}{2}-\frac{1}{2N}\bigg).\log(k)\Bigg)} \\[1ex]
\geq M^{-2}\chi\Big(\{\text{Re}C_{ij}\},\{\text{Im}C_{ij}\}:(i,j)\in
\Phi\Big)
\end{multline*}
Then
\begin{align*}
\mathbb{P}_{\epsilon}(Z_H:D)& \geq
\liminf_{k}{(Mk)^{-2}\log(E_{\epsilon}(f_k))}+\Big (2-\frac{1}{N} \Big)|\log
\epsilon| \\
&\quad+2\iint_{K_N}\log(\max(H(s,t),\epsilon))dsdt + L_3\nonumber
\end{align*}
The eigenvalue distribution of $f_k$ equals that of $z_{k,\epsilon}$
and converges as $k\to+\infty$ to the measure $\sigma_{\epsilon,\frac1{\sqrt{N}}}$, we may apply
Lemma~2.3 of~\cite{DT} concerning the asymptotics of $E_\epsilon(f_k)$ as $k\to\infty$.
Using also Lemma~\ref{prop1}, we get
\[
\delta_0(Z_H:D)= \limsup_{\epsilon\to0}{\frac{\mathbb{P}_\epsilon(Z_H:D)}{|\log \epsilon|}}
\geq 1 +2\cdot\mathrm{area}(\mathrm{supp}(H)\cap K_N).
\]
Taking $N$ arbitrarily large completes the proof.
\end{proof}

The following Theorem gives us an upper bound on $\delta_0(Z_{H}:D)$
without any conditions on the support of $H$.

\begin{teo}\label{main2}
Let $H\geq 0$, $H\in L^{1}([0,1]^2)$ have essentially bounded coordinate expectations $CE_1(H)$
and $CE_2(H)$, as in equations~\eqref{eq:CE}.
Then
\[
\delta_0(Z_{H}:D)\leq
\min\{\,2\,,\,\,1+2\,\mathrm{area}(\mathrm{supp}(H))\}.
\]
\end{teo}

\begin{proof}
First of all it is clear that $\delta_0(Z_{H}:D)\leq
\delta_0(Z_{H})\leq 2$.\\
By standard arguments we can find $\omega$ in regular block form such that both $\norm{Z_{H}-Z_{\omega}}_2$ 
and $\mathrm{area}(\mathrm{supp}(H)\triangle\mathrm{supp}(w))$ are arbitrarily small. Using this, given $\delta>0$ we can find projections
$p_1,\, q_1,\, p_2,\, q_2,\ldots, p_n,\,
q_n$ in $W^{*}(D)$ such that if $i\neq j$, then
$p_{i}\otimes q_{i}$ is orthogonal to $p_{j}\otimes q_{j}$ in
$W^{*}(D)\overline{\otimes} W^{*}(D)$  and such that
\begin{align}\label{tau1}
\sum_{i=1}^{n}{\tau(p_i)\tau(q_i)}&>1-\mathrm{area}(\mathrm{supp}(H))-\delta/3 \\
\sum_{i=1}^{n}{\norm{p_iZ_{H}q_i}_2}&<\delta/4.
\end{align}
Take $R>\max\{\norm{Z_H}_2,\norm{D}_2\}$. Using Lemma 2.9 of \cite{tubularity}, given $\epsilon>0$ there exist $m_0,
\gamma_0, k_0$ such that for $m\geq m_0$, $\gamma<\gamma_0$,
$k\geq k_0$ and for every
$(A,B) \text{\,\,and\,\,} (\tilde{A},\tilde{B})\in\Gamma_{R}(Z_H,D;m,k,\gamma)$ there exists a unitary $U\in M_{k}(\C)$ such that
\begin{equation}\label{approx}
\norm{U\tilde{B}U^{*}-B}_2<\epsilon.
\end{equation}
For $m$ and $k$ sufficiently big and $\gamma$ sufficiently small we can find spectral projections of $B$ 
\[
P_1,Q_1,\ldots,P_n,Q_n\in M_{k}(\C)
\]
and spectral projections of $\tilde{B}$
\[
\tilde{P}_1,\tilde{Q}_1,\ldots,\tilde{P}_n,\tilde{Q}_n\in M_{k}(\C)
\]
such that if $i\neq j$ then $P_{i}\otimes Q_{i}$
is orthogonal to $P_{j}\otimes Q_{j}$ in $M_{k}(\C)\otimes M_{k}(\C)$ and $\tilde{P}_{i}\otimes \tilde{Q}_{i}$ is orthogonal to
$\tilde{P}_{j}\otimes \tilde{Q}_{j}$
satisfying 
\[
|\tr_k(P_i)-\tau(p_i)|<\frac\delta{3n},\quad
|\tr_k(Q_i)-\tau(q_i)|<\frac\delta{3n},\quad
\sum_{i=1}^{n}\norm{P_iAQ_i}_2<\frac\delta2
\]
\[
|\tr_k(\tilde{P}_i)-\tau(p_i)|<\frac\delta{3n},\quad
|\tr_k(\tilde{Q}_i)-\tau(q_i)|<\frac\delta{3n},\quad
\sum_{i=1}^{n}\norm{\tilde{P}_i\tilde{A}\tilde{Q}_i}_2<\frac\delta2.
\]
Taking $\epsilon$ sufficiently small and using (\ref{approx}) together with the fact that we can always approximate 
these projections with polynomials in $B$ and $\tilde{B}$ in the $|\cdot|_2$, we can also guarantee that
\[\norm{P_i-U\tilde{P_i}U^{*}}_2<\frac{\delta}{6nR}, \quad
\norm{Q_i-U\tilde{Q_i}U^{*}}_2<\frac{\delta}{6nR} \quad(1\leq i\leq n).
\]
Therefore, 
\begin{equation}\label{delta-ngbd}
\sum_{i=1}^{n}{\norm{P_i(U\tilde{A}U^{*})Q_{i}}}_2<\sum_{i=1}^{n}{\Bigg(\frac{3\delta\norm{\tilde{A}}}{6nR}+\norm{\tilde{P}_i\tilde{A}\tilde{Q}_i}_2\Bigg)}
<\delta.
\end{equation}
Let $\Omega_{R}(H,k)=\{X\in
M_{k}(\C):\norm{X}_2\leq
R,\,\,\,P_iXQ_i=0\,\,\,\text{for}\,\,i=1,\ldots,n\}$, this is a
ball of radius $R$ in a space of real dimension $d(k)=
2k^2(1-\sum_{i=1}^{n}{\tr_{k}(P_i)\tr_{k}(Q_i)})$.
By (\ref{delta-ngbd}) it is clear that
\begin{equation}\label{Omega}
\Gamma_{R}(Z_H:D\,;m,k,\gamma)\subseteq
\theta(N_{\delta}(\Omega_{R}(H,k)))
\end{equation}
where $\theta(N_{\delta}(\Omega_{R}(H,k)))$ is the unitary orbit
of the $\delta-$neighborhood of $\Omega_{R}(H,k)$.
Taking
the $P_\delta$ packing number on both sides of~\eqref{Omega}, we get
$$P_{\delta}(\Gamma_{R}(Z_H:D\,;m,k,\gamma))\leq
P_{\delta}(\theta(N_{\delta}(\Omega_{R}(H,k))))
\leq P_{\delta}(U_{k}(\C))\cdot
P_{\delta}(N_{\delta}(\Omega_{R}(H,k))).$$ Using Theorem~7 of~\cite{Sz},
there exists a constant $K_1$ independent of $k$ such
that
\begin{equation}
P_{\delta}(U_{k}(\C))\leq \Bigg(\frac{K_1}{\delta}\Bigg)^{k^2}.
\end{equation}
On the other hand, standard packing number estimations gives us
\begin{equation}
P_{\delta}(N_{\delta}(\Omega_{R}(H,k)))\leq
P_{\delta}(\Omega_{R+\delta}(H,k))\leq
\Bigg(\frac{K_2(R+\delta)}{\delta}\Bigg)^{d(k)}
\end{equation}
where $K_2$ is a constant independent of $k$. It follows that
$$P_{\delta}(\Gamma_{R}(Z_H:D\,;m,k,\gamma))\leq
\Bigg(\frac{K_1}{\delta}\Bigg)^{k^2}\cdot
\Bigg(\frac{K_2(R+\delta)}{\delta}\Bigg)^{d(k)}.$$ Now using
(\ref{tau1}) yields
\begin{align*}
\frac{d(k)}{k^2}&=2\Big(1-\sum_{i=1}^{n}{\tr_{k}(P_i)\tr_{k}(Q_i)\Big)}\leq
2\Big(1-\sum_{i=1}^{n}{\tau(p_i)\tau(q_i)}+2\delta/3\Big) \\
&\leq2\Big(\mathrm{area}(\mathrm{supp}(H))+\delta\Big).
\end{align*}
Therefore,
\begin{eqnarray*}
\limsup_{k}\frac{1}{k^2}\log(P_{\delta}(\Gamma_{R}(Z_H:D\,;m,k,\gamma)))
&\leq& \log(K_1)+|\log(\delta)| + \\
&+& 2(\mathrm{area}(\mathrm{supp}(H))+\delta)\cdot\log(K_2(R+\delta))\\
&+& 2(\mathrm{area}(\mathrm{supp}(H))+\delta)\cdot |\log(\delta)|.
\end{eqnarray*}
When $\gamma\to 0$ and $m\to +\infty$ we obtain
$$\mathbb{P}_\delta(Z_H:D)\leq
(1+2\cdot\mathrm{area}(\mathrm{supp}(H))+2\delta)\cdot|\log\delta|+C$$
where $C$ is a constant. It follows that
$$\delta_{0}(Z_H:D)=\limsup_{\delta\to
0}{\frac{\mathbb{P}_\delta(Z_H:D)}{|\log \delta|}}
\leq 1+2\cdot\mathrm{area}(\mathrm{supp}(H)).$$

\end{proof}
\section{Concluding remarks and questions}
\label{sec:concl}

Since the free entropy dimension of $Z_H$ in the presence of $D$ is
a lower bound for the free entropy dimension of $Z_H$,
from Theorems~\ref{main} and~\ref{main2} we have that for any
$H$ as in Theorem~\ref{main},
\begin{equation}\label{eq:pres}
1+2\,\mathrm{area}(\mathrm{supp}(H))=
\delta_{0}(Z_{H}:D)\leq
\delta_{0}(Z_H).
\end{equation}
However,
$1+2\,\mathrm{area}(\mathrm{supp}(H))$ is not the actual value
of $\delta_{0}(Z_H)$ in all cases.
For example, if $n\ge2$ and if
$H$ is the characteristic function of $\cup_{i=1}^nT_i$,
where $T_i=\{(x,y)\in[0,1]:\,\frac{i-1}{n}\leq x<y\leq \frac{i}{n}\}$,
then the moments of $Z_H$ agree with the moments of a nonzero multiple of the
quasinilpotent DT--operator $T$.
Therefore, in
this case we have
\begin{equation}
\delta_{0}(Z_{H}:D)=1+\frac{1}{n}<\delta_{0}(Z_H)=\delta_0(T)=2.
\end{equation}

Of course, if $D$ belongs to the von Neumann algebra generated by $Z_H$, then
equality holds in~\eqref{eq:pres}.
It is an interesting question, when do we have $D\in W^*(\{Z_H\})$?
More generally, what is the von Neumann algebra generated by $Z_{H}$?
When is it a factor?
Is it then
an interpolated free group factor?
A particular case of interest is 
when $H$ is the characteristic function of the band
\[
\{(x,y)\mid 0\le x<y<\min(1,x+\alpha)\},
\]
for $\alpha\in(0,1)$, as is drawn in Figure~\ref{fig:band} (on page~\pageref{fig:band}).

\end{document}